	\documentclass{article}
        \usepackage{latexsym}

	\newcommand{\al}{\alpha}
	
	\newcommand{\eps}{\epsilon}
	\newcommand{\gam}{\gamma}
	\newcommand{\lam}{\lambda}
	\newcommand{\vph}{\varphi}
	\newcommand{\sig}{\sigma}
	\newcommand{\tht}{\theta}
	\newcommand{\om}{\omega}
	\newcommand{\Del}{{\mathit \Delta}}
	\newcommand{\Gam}{{\mathit \Gamma}}
	
	\newcommand{\Sig}{{\mathit \Sigma}}
	\newcommand{\Om}{{\mathit \Omega}}

	\font\bbb=msbm7 
	\font\BBB=msbm10 
	\newcommand{\NN}{{\mbox{\BBB{N}}}}
	\newcommand{\ZZ}{{\mbox{\BBB{Z}}}}

	\newcommand{\re}{{\mbox{\bbb{R}}}}
	\newcommand{\RE}{{\mbox{\BBB{R}}}}
	\newcommand{\co}{{\mbox{\bbb{C}}}}
	\newcommand{\CO}{{\mbox{\BBB{C}}}}
	\font\frak=eufm10 at 11 pt

        \newcommand{\gt}{\mbox{\frak{t}}}
	\newcommand{\OO}{{\cal O}}
	\newcommand{\EE}{{\cal E}}
	\newcommand{\HH}{{\cal H}}
	\newcommand{\medwedge}{\mbox{\fontsize{12pt}{0pt}\selectfont $\wedge$}}

	\newcounter{sect}
	\newcommand{\sect}[1]{\vspace{4ex}\addtocounter{sect}{1}
		\begin{flushleft}
		{{\large\bf \arabic{sect}. {#1}}}
		\end{flushleft}
		\setcounter{thm}{0}
		\setcounter{equation}{0}
		\def\theequation{\arabic{sect}.\arabic{equation}}}
	\newcommand{\subsect}[2]{\vspace{.5ex}
		\begin{flushleft}
		{{\bf \arabic{sect}.{#1} {#2}}}
		\end{flushleft}}
	\newtheorem{thm}{Theorem}[sect]
	\newtheorem{prop}[thm]{Proposition}
	\newtheorem{lemma}[thm]{Lemma}
	\newtheorem{cor}[thm]{Corollary}
	\newtheorem{defn}[thm]{Definition}
	\newtheorem{rmk}[thm]{Remark}
	\newtheorem{ex}[thm]{Example}
	\newcommand{\proof}[1]{\noindent {\em Proof.}$\quad$ {#1} $\hfill\Box$
				\vspace{2ex}}
	\newcommand{\be}{\begin{equation}}
	\newcommand{\ee}{\end{equation}}
	\newcommand{\bea}{\begin{eqnarray}}
	\newcommand{\eea}{\end{eqnarray}}
	\newcommand{\nno}{\nonumber \\}
	\newcommand{\sep}[1]{\!\!\!\! &{#1}& \!\!\!\! }
	\newcommand{\eq}{\sep{=}}
	\newcommand{\vc}{\sep{ }}

	\newcommand{\bra}{\langle}
	\newcommand{\ket}{\rangle}
	\newcommand{\m}{\backslash}
	\newcommand{\inv}[1]{\frac{1}{#1}}
	\newcommand{\hf}{{\textstyle\inv{2}}}
	\newcommand{\e}[1]{e^{{#1}}}
	\newcommand{\ii}{\sqrt{-1}}
	\newcommand{\dr}{d}
	\newcommand{\ch}{\,\mathrm{char}\,}
	\newcommand{\td}{\,\mathrm{td}\,}
	\newcommand{\tr}{\,\mathrm{tr}}
	\newcommand{\supp}{\,\mathrm{supp}}
	\newcommand{\mult}{\,\mathrm{mult}}
	\newcommand{\set}[2]{\{{#1}\,|\,{#2}\}}
	\newcommand{\norm}[1]{||{#1}||}

        \newcommand{\lat}{{\cal L}}
        \newcommand{\ring}{\CO[\lat^*]}
	\newcommand{\kb}{{\bar{k}}}
	\newcommand{\zb}{{\bar{z}}}
	\newcommand{\rk}{_{r,k}}
	\newcommand{\kl}{^{k,l}}
	\newcommand{\lrk}{{\lam\rk}}
	\newcommand{\ncr}{{\nu^C_r}}
	\newcommand{\Ncr}{N^{C,(1,0)}_r}
	\newcommand{\Nmcr}{N^{-C,(1,0)}_r}
	\newcommand{\ECr}{K^C(N_r)\otimes E|_{F_r}}
	\newcommand{\Epr}{K^+(N_r)\otimes E|_{F_r}}
	\newcommand{\Emr}{K^-(N_r)\otimes E|_{F_r}}
	\newcommand{\txi}{_{T,\xi}}
	\newcommand{\sumi}{\sum_{i=1}^{2n}}
	
	\newcommand{\sumk}{\sum^n_{k=0}}
	\newcommand{\sumkr}{\sum^{n_r}_{k=0}}
	\newcommand{\sumr}{\sum^m_{r=1}}
	
	\newcommand{\two}[4]{\left\{    \begin{array}{ll}
					{#1}, & {\mbox{if }} {#2}, \\
					{#3}, & {\mbox{if }} {#4}
					\end{array}     \right.}
	\newcommand{\ka}{K\"ahler }
	\newcommand{\hme}{{H^k(M,\OO(E))}}
	\newcommand{\pdr}{\partial}
	\newcommand{\pdrb}{\bar{\pdr}}
	\newcommand{\tina}{\tilde{\nabla}}
	\newcommand{\ha}{^{1,0}}
	\newcommand{\ah}{^{0,1}}
	\newcommand{\cht}{\mathrm{ch}_T}
	\newcommand{\dol}{\medwedge(T^{*(0,1)}M)\otimes E}
	\newcommand{\bfh}{{\mbox{{\bf H}}}}
	\newcommand{\tilpi}{\tilde{\pi}}
	\newcommand{\dv}{\dr v}
	\newcommand{\grad}{\mathrm{grad}}
	\newcommand{\spec}{\mathrm{Spec}}

\begin{document}
$\!\,{}$

        \vspace{-5ex}

        \begin{flushright}
{\tt dg-ga/9701009} (January, 1997)\\
Adelaide IGA preprint 1997-01
        \end{flushright}

	\begin{center}
	{\LARGE\bf Equivariant Holomorphic Morse Inequalities III:\\
	\vspace{1ex}
		Non-Isolated Fixed Points}\\

	\vspace{4ex}
	{\large\rm Siye Wu}\footnote{Present Address: School of Mathematics,
	Institute for Advanced Study, Olden Lane, Princeton, NJ 08540, USA.
	E-mail addresses: {\tt swu@maths.adelaide.edu.au, swu@math.ias.edu}}\\
	{\em Department of Pure Mathematics, The University of Adelaide,
	Adelaide, South Australia 5005, Australia}

	\vspace{2ex}

	{\large\rm Weiping Zhang}\\
	{\em Nankai Institute of Mathematics, Tianjin, 300071, P.\ R.\ China}

	\end{center}

	\vspace{3ex}

	\begin{quote}
{\small {\bf Abstract.} 
We prove the equivariant holomorphic Morse inequalities for a holomorphic 
torus action on a holomorphic vector bundle over a compact \ka manifold
when the fixed-point set is not necessarily discrete.
Such inequalities bound the twisted Dolbeault cohomologies of the \ka 
manifold in terms of those of the fixed-point set.
We apply the inequalities to obtain relations of Hodge numbers of the
connected components of the fixed-point set and the whole manifold.
We also investigate the consequences in geometric quantization, especially
in the context of symplectic cutting.}

	\end{quote}

	\vspace{2ex}

\sect{Introduction}
In [\ref{MW}] a rigorous proof of equivariant holomorphic Morse inequalities
of Witten [\ref{W}] was found based on the heat kernel approach.
The result proved holds for circle actions on holomorphic vector bundles over
compact \ka manifolds with isolated fixed points.
In [\ref{WuII}], this result is extended to torus and non-Abelian group
actions and it was shown that the \ka condition is necessary.
The purpose of this paper is to prove the general case when the fixed points
can be non-isolated.

Recall that in his heat kernel proof of the Bott-Morse inequalities 
for the de Rham complex, Bismut [\ref{Bi}] used a further rescaling
near the critical set of a Morse function in the sense of Bott,
besides the deformation of Witten [\ref{W0}].
In our treatment of the holomorphic analog with non-isolated fixed points,
instead of using [\ref{Bi}], we adapt the methods and techniques in
[\ref{BL}], where a general and direct localization procedure is developed
which applies to a wide range of problems on localization in index theory.
For example, to prove the Bott-Morse inequalities using this method, 
no rescaling near the critical set is needed.

As in [\ref{WuII}], the equivariant holomorphic Morse inequalities for
circle actions imply those for torus actions;
we derive the latter by considering various circle subgroups of the torus.
For circle actions, we deform the Dolbeault operator as in [\ref{W}].
The corresponding Laplacian is roughly of the form $\Del+|V|^2+L_V$,
where $\Del$ is the standard Laplacian, $V$ is the vector field that generates
the circle action, and $L_V$, the Lie derivative.
$L_V$ is a first order differential operator that commutes with all other
operators..
The crucial step is to restrict the problem to the eigenspaces of $L_V$,
on each of which $L_V$ is a constant and therefore the techniques of [\ref{BL}]
can be applied.
This leads naturally to an equivariant Morse theory.

The equivariant holomorphic Morse inequalities relate the Dolbeault
cohomologies (twisted by a holomorphic vector bundle) of the whole manifold
to those of the fixed point set. 
We apply our result to two situations.
First, considering the exterior power of the holomorphic tangent bundle,
we obtain relations of Hodge numbers of the fixed submanifolds, including
some results of Carrell-Sommese [\ref{CS}] and Carrell-Lieberman [\ref{CL}].  
Second, we apply the inequalities to symplectic cutting and geometric
quantization.
We obtain a gluing formula for the Poincar\'e-Hodge polynomials under 
the symplectic cutting precess of Lerman [\ref{L}].
Applying our result to the pre-quantum line bundles, we recover and/or 
strengthen a few results in [\ref{DGMW}, \ref{MS}, \ref{TZ}] 
when the symplectic manifold is K\"ahler. 

Throughout this paper, $\NN$, $\RE$, $\RE^\pm$ and $\CO$ denote the sets of 
non-negative integers, real numbers, positive (negative) real numbers and
complex numbers, respectively.

\sect{Main Results}
In this section, we state the equivariant holomorphic Morse inequalities
for holomorphic torus action on holomorphic vector bundle over a 
compact \ka manifold when the fixed-point set need not be isolated.
We then show that the result for torus action can be deduced from that of
circle action.
The latter will be proved in Section 3.

\subsect{1.}{Equivariant holomorphic Morse inequalities for torus actions}
We first recall from [\ref{WuII}] a few notations of the Lie algebra and
the formal character ring of the torus group.
Let $\gt=\ii\,\mathrm{Lie}(T)$ be the Lie algebra of the torus group $T$
and $\lat$, the integral lattice in $\gt$;
the dual lattice $\lat^*$ in $\gt^*$ is the weight lattice.

\begin{defn}
Let $\ring$ be the formal character ring of $T$ consisting of elements 
$q=\sum_{\xi\in\lat^*}q_\xi\e{\xi}$ ($q_\xi\in\ZZ$).
We say $q\ge0$ if $q_\xi\ge0$ for all $\xi\in\lat^*$.
Let $Q(t)=\sumk q_kt^k\in\ring[t]$ be a polynomial of degree $n$ with
coefficients in $\ring$, we say $Q(t)\ge0$ if $q_k\ge0$ in $\ring$ for all $k$.
For two such polynomials $P(t)$ and $Q(t)$, we say $P(t)\le Q(t)$ if
$Q(t)-P(t)\ge0$.
\end{defn}

If $W$ is a finite dimensional representation of $T$, We denote by 
$\mult_\xi(W)$ the multiplicity of the weight $\xi\in\lat^*$ in $W$.
The character of $W$ is $\ch(W)=\sum_{\xi\in\lat^*}\!\mult_\xi(W)\,\e{\xi}\ge0$
in $\ring$.
Let the {\em support} of $W$ be
$\supp(W)=\set{\xi\in\lat^*}{\mult_\xi(W)\ne 0}$.
For any $\tht\in\gt^*$, there is a homomorphism $\ring\to\CO$
given by $\e{\xi}\mapsto\e{\ii\bra\xi,\tht\ket}$.
For instance, $\ch(W)\mapsto\tr_W\e{\ii\tht}$ under this homomorphism.

Now let $(M,\om)$ be a compact \ka manifold of complex dimension $n$.
Suppose that $T$ acts on $M$ effectively preserving the complex structure
and the \ka form.
If the fixed-point set $F\subset M$ of the $T$-action is non-empty,
then the action is Hamiltonian [\ref{Fr}].
$F$ is a finite union of connected compact \ka submanifolds $F_1,\dots,F_m$;
let $n_1,\dots,n_m$ be their complex dimensions, respectively.
For each connected component $F_r$ ($1\le r\le m$), the complexification
of the normal bundle $N_r\to F_r$ in $M$ has the decomposition
$N^\co_r=N\ha_r\oplus N\ah_r$, and $N\ha_r$ is a holomorphic vector bundle
over $F_r$ of rank $n-n_r$.
The torus $T$ acts on $N_r$ preserving the base points in $F_r$.
Moreover, the weights of the isotropy representation on the normal fiber
remains unchanged within any connected component of $F$.
Let $\lam_{r,k}$ ($1\le k\le n-n_r$) be the isotropy weights on $N_r$.
$N_r$ splits, and $N\ha_r$ splits holomorphically, into the direct
sum of various sub-bundles, each with a given weight.
The hyperplanes $(\lam\rk)^\perp\subset\gt$ cut $\gt$ into open polyhedral
cones called {\em action chambers}, as in [\ref{GLS},\ref{PW},\ref{WuII}].
We fix a {\em positive} action chamber $C$.
Let $\lam^C\rk=\pm\lrk$ be the {\em polarized weights},
with the sign chosen so that $\lam^C\rk\in C^*$. 
(Here $C^*$ is the dual cone in $\gt^*$ defined by
$C^*=\set{\xi\in\gt^*}{\bra\xi,C\ket>0}$.)
We define the {\em polarizing index} $\ncr$ of the component $F_r$
with respect to $C$ as the number of weights $\lrk\in-C^*$.
Let $N^C_r$ be the direct sum of the sub-bundles corresponding to
the weights $\lrk\in C^*$.
$\ncr$ is the rank of the holomorphic vector bundle $N^{-C,(1,0)}_r$;
that of $N^{C,(1,0)}_r$ is $\nu^{-C}_r=n-n_r-\ncr$.
Finally, we define the {\em polarized symmetric tensor product}
(with respect to $C$) of the holomorphic normal bundle $N\ha_r$ by
	\be
K^C(N_r)=S((\Ncr)^*)\otimes S(\Nmcr)\otimes\medwedge^{\ncr}(\Nmcr).
	\ee
This is bundle of eigenspaces of small eigenvalues implicitly contained in
[\ref{W}].
It also appeared in [\ref{CG}], which generalizes [\ref{GLS}] to cases
with non-isolated fixed points.

If $E$ is a holomorphic vector bundle over $M$ on which the $T$-action
lifts holomorphically, then the fiber $E_p$ over each fixed point $p\in F$ 
is a representation of $T$, and $\ch(E_p)$ is also constant within
any connected component of $F$.
Consider an infinite dimensional bundle holomorphic $\ECr$.
The torus $T$ acts on the total space while preserving the base points
in $F_r$.
It is easy to see that each sub-bundle of any given weight is 
a holomorphic vector bundle of finite rank, i.e.,
$\ECr=\oplus_{\xi\in\lat}E^C_{r,\xi}$, where $E^C_{r,\xi}$ is a $T$-invariant 
sub-bundle of finite rank on which the torus acts with weight $\xi$.
The cohomology groups $H^k(F_r,\OO(\ECr))$ are the sum of those with
coefficients in $E^C_{r,\xi}$, each equipped an induced $T$-action.
Therefore $\ch H^k(F_r,\OO(\ECr))=\sum_{\xi\in\lat}\dim_\co 
H^k(F_r,\OO(E^C_{r,\xi}))\,\e{\xi}$ is a well-defined element in $\ring$.
Moreover, $\supp H^k(F_r,\OO(\ECr))$ is contained in a suitably shifted
cone $-C^*$ in $\gt^*$.

Our main result is the following

\begin{thm}\label{TORUS}
For each choice of the positive action chamber $C$, we have
the strong equivariant holomorphic Morse inequalities
        \be\label{strong}
\sumr t^{\ncr}\sumkr t^k\ch H^k(F_r,\OO(\ECr))
=\sumk t^k\ch H^k(M,\OO(E))+(1+t)Q^C(t)
        \ee
for some $Q^C(t)\ge0$ in $\ring[t]$.
\end{thm}

\begin{rmk} {\em
1. Formula (\ref{strong}) clearly implies the corresponding weak inequalities
        \be\label{weak}
\ch H^k(M,\OO(E))\le\sumr\ch H^{k-\ncr}(F_r,\OO(\ECr)).
        \ee
2. It is easy to see that for any choice of $C$, (\ref{strong}) reduces to
the Atiyah-Bott-Segal-Singer fixed-point theorem [\ref{AB}, \ref{ASS}]
after setting $t=-1$.
In fact,
        \bea
\sumk(-1)^k\ch\hme
\eq\sumr\int_{F_r}\cht(\ECr)\td(F_r)				\nno
\eq\sumr\int_{F_r}\cht\left(\frac{E|_{F_r}}{\det(1-N^{*(1,0)}_r)}\right)
	\td(F_r),				\label{abss}
        \eea
where $\cht$, $\td$ stand for the equivariant Chern character 
and the Todd class, respectively.\\
3. If $p\in F$ is an isolated fixed point, then the normal bundle
$N_p=T_pM$ and
	\be
\ch H^0(\{p\},\OO(K^C(N_p)\otimes E_p))=\ch(E_p)
\prod_{\lam^C_{p,k}\in C^*}\inv{1-\e{-\lam_{p,k}}}
\prod_{\lam^C_{p,k}\in-C^*}\frac{\e{-\lam^C_{p,k}}}{1-\e{-\lam^C_{p,k}}}.
	\ee
Therefore, (\ref{strong}) reduces to the result in [\ref{WuII}] when
$F$ is discrete.\\
4. Though the \ka assumption is not necessary in the fixed-point formula
(\ref{abss}), it is essential for the strong inequalities (\ref{strong}) even
when all the fixed points are isolated and when $T$ is a circle group
[\ref{WuII}].}
\end{rmk}

As in [\ref{WuII}], (\ref{weak}) or (\ref{strong}) 
is a set of inequalities for each choice of the action chamber $C$.
These inequalities provides bounds, along various directions in $\gt^*$
given by $C$, the multiplicities of weights in the Dolbeault cohomology
group $H^k(M,\OO(E))$ ($0\le k\le n$) in terms of the fixed-point data,
which includes $F_r$ and the bundles $N_r, E|_{F_r}\to F_r$ with $T$-actions.
The applications will be given in Section 4.

When the group acting on $M$ is a non-Abelian group $G$, we can apply
(\ref{strong}) to the maximal torus $T$ of $G$.
There is in addition an action of the Weyl group $W$ on the fixed-point set $F$
of $T$.
Each $w\in W$ induces an action on the set of connected components
$\{F_1,\dots,F_m\}$ of $F$.
The sum over the connected components on the left hand side of (\ref{strong})
can be rearranged into sums over $W$ (after incorporating the character of 
the isotropy representation of $T$ on the bundle $E|_F$)
and over its orbits in the set of connected components.
Thus we can obtain the non-Abelian version of the holomorphic Morse
inequalities like the case when $F$ is isolated [\ref{WuII}].
The details are left to the interested reader.

\subsect{2.}{Reduction to the case of circle actions}
When the torus $T$ is the circle group $S^1$,
the Lie algebra is $\mathrm{Lie}(S^1)=\ii\RE$.
There are two action chambers $\RE^\pm$, labeled by $\pm$ for simplicity.
The weights of isotropy representation of $S^1$ on $N_r$ ($1\le r\le m$)
are $\lam_{r,k}\in\ZZ\m\{0\}$ ($1\le k\le n_r$).
As before, $N_r^\pm$ is the direct sum of the sub-bundles
corresponding to the positive and negative weights, respectively,
and $K^\pm(N_r)=S((N^{\pm,(1,0)}_r)^*)\otimes S(N^{\mp,(1,0)}_r)\otimes
\medwedge^{\mathrm{top}}(N^{\mp,(1,0)}_r)\otimes E|_{F_r}$.
Let $\nu_r$ be the rank of $N^{-,(1,0)}_r$;
that of $N^{+,(1,0)}_r$ is then $n-n_r-\nu_r$.

Then we have the following result for the $S^1$-case.

\begin{thm}\label{CIRCLE}
Under the above assumptions, we have the 
strong equivariant holomorphic Morse inequalities for circle actions
        \be\label{strong1+}
\sumr t^{\nu_r}\sumkr t^k\ch H^k(F_r,\OO(\Epr))
=\sumk t^k\ch H^k(M,\OO(E))+(1+t)Q^+(t),
	\ee
        \be\label{strong1-}
\sumr t^{n-n_r-\nu_r}\sumkr t^k\ch H^k(F_r,\OO(\Emr))
=\sumk t^k\ch H^k(M,\OO(E))+(1+t)Q^-(t),
	\ee
where $Q^\pm(t)\ge0$.
\end{thm}

Clearly (\ref{strong1+}) follows from (\ref{strong1-})
by reversing the $S^1$-action.
Section 3 will be devoted to proving (\ref{strong1-}).
To show that the strong inequalities (\ref{strong}) of the torus case follows
from (\ref{strong1+}) or (\ref{strong1-}), we proceed slightly differently
from [\ref{WuII}].

Recall that a cone $C$ in $\gt$ is {\em proper} if there is $\xi\in\gt^*$
such that $\bra\xi,C\m\{0\}\ket>0$.
If $C$ is an open proper cone, so is $C^*$ in $\gt^*$.
For example, if $T$ acts on $M$ effectively, then all action chambers
are open proper cones.
(See [\ref{PW}] for the geometry of cones in the context of Hamiltonian
torus actions.)

\begin{lemma}
If $C$ be an open proper cone in $\gt$, then for any $\xi\in C^*\cap\lat^*$,
there is an element $v\in C\cap\lat$
such that the hyperplane $\xi+v^\perp=\set{\lam\in\gt^*}{\bra\lam-\xi,v\ket=0}$
contains no points in $C^*\cap\lat^*$ other than $\xi$.
\end{lemma}

\proof{Choose an open proper subcone $D\subset C$ such that
$\overline{D}\subset C\cup\{0\}$.
If $v\in D\cap\lat$, then the intersection $(\xi+v^\perp)\cap(C^*\cap\lat^*)$
is contained in $(C^*\m(\xi+\overline{D^*}))\cap\lat^*$;
the latter is a finite set since $\overline{C^*}\subset D^*\cup\{0\}$.
For each $\lam\in(C^*\m(\xi+\overline{D^*}))\cap\lat^*$, consider
the hyperplane $H_\lam=\set{v\in\gt}{\bra\lam-\xi,v\ket=0}$ in $\gt$.
Let $\pi\colon\gt\m\{0\}\to P(\gt)$ be the canonical projection to
the real projective space $P(\gt)$.
Since the images $\pi(H_\lam)\subset P(\gt)$ are of codimension $1$,
$\pi(\lat)$ is dense in $P(V)$, and $\pi(D)$ is open in $P(\gt)$,
there is an element $v\in\gt$ such that $\pi(v)\in\pi(D\cap\lat)\m\pi(H_\lam)$
for any $\lam\in(C^*\m(\mu+\overline{D^*}))\cap\lat^*$.
By multiplying with a non-zero real number, we may choose $v\in D\cap\lat$.
The condition $n\not\in H_\lam$ means that $\lam+v^\perp$ does not contain
the point $\lam$.
The result follows.}

\noindent {\em Reduction of Theorem~\ref{TORUS} to Theorem~\ref{CIRCLE}:}
$\quad$ 
For any action chamber $C$, there is $\lam_0\in\lat^*$ such that
	\be
\supp\,H^k(F_r,\OO(\ECr)), \supp\,\hme\subset(\lam_0-C^*)\cap\lat^*
	\ee
for all $1\le r\le m$, $0\le k\le n$.
(\ref{strong}) is equivalent to
	\be\label{mult}
\sumr t^{\ncr}\sumkr t^k\mult_\xi H^k(F_r,\OO(\ECr))
=\sumk t^k\mult_\xi\hme+(1+t)Q^C_\xi(t),
	\ee
where $Q^C_\xi(t)\ge0$ for all $\xi\in(\lam_0-C^*)\cap\lat^*$.
For any $v\in C\cap\lat$, let $S^1$ be the circle subgroup generated by $v$.
Consider this circle action on $M$ and $E$.
The isotropy weights on $F_r$ ($1\le r\le m$) are $\bra\lam_{r,k},v\ket$
($1\le k\le n_r$).
Since $\bra\lam^C_{r,k},v\ket>0$, $K^C(N_r)=K^+(N_r)$ for all $1\le r\le m$.
By the above Lemma, for any $\xi\in(\lam_0-C^*)\cap\lat^*$ there is an element
$v\in C\cap\lat$ such that $(\xi+v^\perp)\cap((\lam_0-C^*)\cap\lat^*)=\{\xi\}$.
It follows that if $W$ is a representation of $T$ such that 
$\supp(W)\subset(\lam_0-C^*)\cap\lat^*$ and the multiplicities of weights
are finite, then $\mult_\xi(W)=\mult_{\bra\xi,v\ket}(W)$.
Therefore (\ref{mult}) follows from (\ref{strong1+}). $\hfill\Box$

\sect{The case of holomorphic circle actions}
The entire section is devoted to proving (\ref{strong1-})
in Theorem~\ref{CIRCLE}.
In subsection 3.1, we study Witten's deformation of the Dolbeault operator
on a compact \ka manifold as well as on flat spaces.
In subsection 3.2, we establish the Taylor expansion of the deformed operator 
near the fixed-point set.
In subsection 3.3, we split the space of sections into a subspace
corresponding to small eigenvalues and its complement, and estimate
the deformed operator and its resolvent on these subspaces.
It is crucial that these estimates are done on each eigenspace of $L_V$.
Finally, in subsection 3.4, we prove Theorem~\ref{CIRCLE}.

\subsect{1.}{Witten's deformation of the Dolbeault operator}
We consider a holomorphic $S^1$-action on a \ka manifold $M$ which preserves
the \ka structure.
Let $\om$, $g^{TM}$ and $J$ be the \ka form, the \ka metric and the complex
structure on $M$, respectively.
We assume that the fixed point set $F$ is non-empty.
In this case, there is a moment map $h\colon M\to\RE$ satisfying
$i_V\om=\dr h$, where $V$ is the vector field on $M$ that generates
the $S^1$-action [\ref{Fr}].
We further assume that the $S^1$-action can be lifted holomorphically
to a holomorphic vector bundle $E$ over $M$.
We can choose an $S^1$-invariant Hermitian form on $E$.
Then the Hermitian connection $\nabla^E$ is also $S^1$-invariant.
$\nabla^E$ induces an ($S^1$-equivariant) action, also denoted by $\nabla^E$,
on the space of $E$-valued differential forms $\Om^*(M,E)$, 
The group elements of $S^1$ acts on the sections of $E$.
Let $L_V$ be the Lie derivative of $E$-valued forms along $V$.
(The fibers of $E$ over different points on the integral curve of $V$
are related by the lifted $S^1$-action.)
Then $-L_V$ is the infinitesimal generator of the $S^1$-action on $\Om^*(M,E)$.
Let $\overline{L}_V=\{i_V,\nabla^E\}$.
Then the operator
	\be\label{rv}
r_V=\overline{L}_V-L_V
	\ee
is an element of $\Gam(M,{\rm End}(E))$.
Over the fiber of a fixed point $p\in F$, $r_V(p)$ is simply
the representation of ${\rm Lie}(S^1)$ on $E_p$;
this is independent of the choice of the connections on $E$.

Let $\pdrb=\pdrb^E$ be the twisted Dolbeault operator of the complex
$\Om^{0,*}(M,E)=\Gam(\medwedge^*(T^{*(0,1)}M)\otimes E)$ and $\pdrb^*$,
its formal adjoint.
Consider Witten's deformation of the (twisted) Dolbeault operator
	\be\label{deform}
\pdrb_h=\e{-h}\pdrb\e{h},	\quad	\pdrb^*_h=\e{h}\pdrb^*\e{-h}.
	\ee
We define the Clifford action.
For $X\in\Gam(T^\co M)$, let $X=X\ha+X\ah$ such that $X\ha\in\Gam(T\ha M)$
and $X\ah\in\Gam(T\ah M)$.
Set 
	\be\label{clifford}
c(X\ha)=\sqrt{2}(X\ha)^*,\quad\quad c(X\ah)=-\sqrt{2}\,i_{X\ah},
	\ee
where $(X\ha)^*\in\Gam(T^{*(0,1)}M)$ corresponds to $X\ha$ via the metric
$g^{TM}$.
It is easy to verify that for $X,Y\in\Gam(T^\co M)$, the anti-commutation
relation $\{c(X),c(Y)\}=-2g^{TM}(X,Y)$.
Let $\{e_i,i=1,\dots,2n\}$ be a (local) orthonormal frame and
$D^M=\sumi c(e_i)\nabla_{e_i}$, the spin$^\co$-Dirac operator
acting on $\Om^{0,*}(M,E)$.

\begin{lemma}
1.	\be
D^M=\sqrt{2}(\pdrb+\pdrb^*);
	\ee
2.	\be
D^M+\ii\,c(V)=\sqrt{2}(\pdrb_h+\pdrb^*_h);
	\ee
3.	\be\label{laplace}
(D^M+\ii\,c(V))^2=(D^M)^2+|V|^2-2\ii L_V
+\hf\ii\sumi c(e_i)c(\nabla_{e_i}V)+\ii\,\tr\,\nabla.V|_{T\ah M}-2\ii\,r_V.
	\ee
\end{lemma}

\proof{Parts 1 and 2 follow from
	\be
\pdrb=\inv{2\sqrt{2}}\sumi c(e_i-\ii Je_i)\nabla_{e_i},\quad
\pdrb^*=\inv{2\sqrt{2}}\sumi c(e_i+\ii Je_i)\nabla_{e_i},
	\ee
	\be
\pdrb_h=\pdrb+\inv{2\sqrt{2}}\sumi c(e_i-\ii Je_i)h_{,i},\quad
\pdrb^*_h=\pdrb^*-\inv{2\sqrt{2}}\sumi c(e_i+\ii Je_i)h_{,i},
	\ee
and that $J\,\grad h=-V$.
To show part 3, we compute
	\be\label{square}
(D^M+\ii\,c(V))^2=(D^M)^2+\ii\{D^M,c(V)\}+|V|^2
	\ee
and
	\be\label{braket}
\{D^M,c(V)\}=\sumi(\{c(e_i),c(V)\}\nabla_{e_i}+c(e_i)c(\nabla_{e_i}V))
=-2\nabla_V+\sumi c(e_i)c(\nabla_{e_i}V).
	\ee
Further, since $L_V=\nabla_V-\nabla.V$ on vector fields,
the induced relation on $\Om^{0,*}(M,E)$ is
	\bea
\nabla_V\eq\overline{L}_V+\inv{4}\sum_{i,j=1}^{2n}g(\nabla_{e_i}\!V,e_j)
	c(e_i)c(e_j)+\inv{2}\tr\,\nabla.V|_{T^{*(0,1)}M}	\nno
	\eq\overline{L}_V+\inv{4}\sumi c(e_i)c(\nabla_{e_i}V)
	-\inv{2}\tr\,\nabla.V|_{T\ah M}.			\label{induce}
	\eea
(\ref{laplace}) follows from (\ref{square}), (\ref{braket}), (\ref{induce})
and (\ref{rv}).}

(\ref{laplace}) was derived in [\ref{MW}] without using the Clifford algebra.

Next we study Witten's deformation on flat spaces [\ref{W}].
Let $W$ be a complex vector space of (complex) dimension $n$ with an Hermitian
form.
Let $\rho$ be a unitary representation of the circle group $S^1$ on $W$
such that all the weights are non-zero.
Suppose $W^\pm$ are the subspaces of $W$ corresponding to positive (negative)
weights, respectively, and $\dim_\co W^-=\nu$, $\dim_\co W^+=n-\nu$.
Let $z=\{z^k\}$ be the complex linear coordinates on $W$ such that 
the Hermitian structure on $W$ takes the standard form and such that
$\rho$ is diagonal with weights $\lam_k\in\ZZ\m\{0\}$ ($1\le k\le n$).
The Lie algebra action is given by the vector field $V=\ii\sum_{k=1}^n
\lam_k(z^k\frac{\pdr}{\pdr z^k}-\zb^k\frac{\pdr}{\pdr\zb^k})$ on $W$.
With respect to the standard \ka form 
$\om=\frac{\ii}{2}\sum_{k=1}^n\dr z^k\wedge\dr\zb^k$, $\rho$ is 
a Hamiltonian action with moment map $h(z)=-\hf\sum_{k=1}^n\lam_k|z^k|^2$.
As in the compact case, we have the deformed operators $\pdrb_h$ and
$\pdrb^*_h$ in (\ref{deform}) and $D+\ii c(V)=\sqrt{2}(\pdrb_h+\pdrb^*_h)$,
where $D$ is the Dirac operator with the same Clifford action (\ref{clifford}).
We set 
$K^\pm(W)=S((W^\pm)^*)\otimes S(W^\mp)\otimes\medwedge^{\mathrm{top}}(W^\mp)$.
Let $E$ be a finite dimensional complex vector space with a Hermitian form
and suppose $E$ carries a unitary representation of $S^1$.

\begin{prop}\label{FLAT}
1. The space of $L^2$-solutions of a given weight of $D+\ii\,c(V)$ on the 
space of $(0,*)$-forms on $W$ with values in $E$ is finite dimensional.
The direct sum of these weight spaces is isomorphic to $K^-(W)\otimes E$ as 
representations of $S^1$.\\
2. When restricted to an eigenspace of $L_V$, the operator $D+\ii\,c(V)$ has
discrete eigenvalues.
\end{prop}

\proof{It suffices to prove the case when $E\cong\CO$ is a trivial
representation.\\
1. Since $W$ decomposes into a direct sum of $1$-dimensional 
representations and since the solution space on the direct sum corresponds
to a tensor product, it suffices to prove the result when $\dim_\co W=1$
with weight $\lam\in\ZZ\m\{0\}$.
Let $\vph\in\Om^{0,*}(\CO)$ solve $(D+\ii c(V))\vph=0$, i.e., 
$\pdrb_h\vph=\pdrb^*_h\vph=0$.
For $\lam<0$, an $L^2$-solution of a fixed weight are proportional to
	\be
\vph_k(z)=\sqrt{\frac{(-\lam)^{k+1}}{\pi k!}}z^k\e{\frac{\lam}{2}|z|^2}
\quad(k\in\NN).
	\ee
$\vph_k$ has unit norm and weight $-k\lam$.
Since $z\in W^*=(W^-)^*$, the direct sum of $\CO\vph_k$ ($k\in\NN$)
is $S((W^-)^*)$.
If $\lam>0$, then $W=W^+$, and a solution of a given weight is proportional to
	\be
\vph_\kb=\sqrt{\frac{\lam^k}{\pi k!}}\,\zb^k\e{-\frac{\lam}{2}|z|^2}\dr\zb
\quad(k\in\NN),
	\ee
of unit norm and weight $(k+1)\lam$.
The direct sum of $\CO\vph_\kb$ ($k\in\NN$) is
$S(\overline{W}^*)\otimes\overline{W}^*\cong S(W^+)\otimes W^+$;
the last isomorphism is induced by the Hermitian form on $W\cong\CO$.\\
2. This follows from (\ref{laplace}) after setting $L_V$ to a constant and
from the standard properties of harmonic oscillators.}

\subsect{2.}{A Taylor expansion of the operator near the fixed-point set}
Since the fixed point set $F$ is the zero set of the Killing vector field $V$,
$F$ is a totally geodesic compact \ka submanifold of $M$.
We denote by $n_F$ and $\nu_F$ the complex dimension and the polarizing index
of $F$, respectively, which are locally constant.
Let $g^{TF}$ be the induced metric and $\dv_F$, the volume element on $F$.
Let $\tilpi\colon N\to F$ be the normal bundle of $F$ in $M$.
We identify $N$ as the orthogonal complement of $TF$ in $TM|_F$,
i.e., $TM|_F=TF\oplus N$ and $g^{TM}=g^{TF}\oplus g^N$,
where $g^N$ is the induced inner product on $N$.
Following [\ref{BL}, Section 8e)], we now describe a coordinate system on 
$M$ near $F$.

If $y\in F$, $Z\in N_y$, let $t\in\RE\mapsto x_t=\exp^M_y(tZ)\in M$ be 
the geodesic in $M$ with $x_0=y$, $\frac{\dr x}{\dr t}|_{t=0}=Z$.
For $\eps>0$, set $B_\eps=\{Z\in N; |Z|<\eps\}$.
Since $M$ and $F$ are compact, there exists $\eps_0>0$ such that for
$0\le\eps<\eps_0$, the map $(y,Z)\in N\mapsto\exp^M_y(Z)\in M$ is a
diffeomorphism from $B_\eps$ onto a tubular neighborhood $U_\eps$ of
$F$ in $M$.						From now on,
we identify $B_\eps$ with $U_\eps$ and use the notation
$x=(y,Z)$ instead of $x=\exp^M_y(Z)$.
Finally, we identify $y\in F$ with $(y,0)\in N$.

Let $\dv_N$ denote the volume element of the fibres in $N$.
Then $\dr v_F(y)\dr v_{N_y}(Z)$ is a natural volume element on
the total space of $N$.
Let $k(y,Z)$ be the smooth positive function defined on $B_{\eps_0}$ 
by the equation $\dv_M(y,Z)=k(y,Z)\dv_F(y)\dv_{N_y}(Z)$.
The function $k$ has a positive lower bound on $U_{\eps_0/2}$.
Also, $k(y)=1$ and $\frac{\pdr k}{\pdr Z}(y)=0$ for $y\in F$;
the latter follows from [\ref{BL}, Proposition 8.9] and the fact 
that $F$ is totally geodesic in $M$.

As in [\ref{BL}, Section 8g)], for $x=(y,Z)\in U_{\eps_0}$, we will identify
$E_x$ with $E_y$ and $\medwedge(T^{*(0,1)}_xM)$ with $\medwedge(T^{*(0,1)}_yM)$
by the parallel transport with respect to the $S^1$-invariant connections 
$\nabla^E$ and $\nabla^{TM}$, respectively, along the geodesic
$t\mapsto(y,tZ)$.
The induced identification of $(\dol)_x$ with $(\dol)_y$ preserves the metric,
the $\ZZ$-grading of the Dolbeault complex, and is moreover $S^1$-equivariant.
Consequently, $D^M$ can be considered as an operator acting on the sections
of the bundle on $\tilpi^*((\dol)|_F)$ over $B_{\eps_0}$.
It still commutes with the $S^1$-action.

For $\eps>0$, let $\bfh(B_\eps)$ (resp.\ $\bfh(N)$) be the set of smooth 
sections of $\tilpi^*(\dol)$ on $B_\eps$ (resp.\ on the total space of $N$).
If $f,g\in\bfh(N)$ have compact support, set
	\be\label{inner}
\bra f,g\ket=\left(\inv{2\pi}\right)^n\int_F
\left(\int_N\bra f,g\ket(y,Z)\dv_{N_y}(Z)\right)\dv_F(y).
	\ee
Notice that the identification of elements in $\bfh(N)$ which have compact
support in $B_{\eps_0}$ with those in $\Gam(\dol)$ with support in $U_{\eps_0}$
is not unitary with respect to the Hermitian product (\ref{inner}).
Consequently $D^M$ as an operator on the sections of $\tilpi^*((\dol)|_F)$
over $B_{\eps_0}$ is not in general self-adjoint with respect to (\ref{inner}).
However $k^\inv{2}D^M k^{-\inv{2}}$ does act as a (formal) self-adjoint
operator on $\bfh(B_{\eps_0})$.

The connection $\nabla^N$ on $N$ induces a splitting $TN=N\oplus T^HN$,
where $T^HN$ is the horizontal part of $TN$ with respect to $\nabla^N$.
Moreover, since $F$ is totally geodesic, this splitting, when restricted
to $F$, is preserved by the connection $\nabla^{TM|_F}$ on $TM|_F$.
If we denote by $p^{TF}$, $p^N$ the orthogonal projections from
$TM|_F$ to $TF$, $N$, respectively, then
$\nabla^{TF}=p^{TF}\nabla^{TM|_F}$ and $\nabla^N=p^N\nabla^{TM|_F}$. 
Let $\tina^F$ be the connection on $(\dol)|_F$ induced by the restrictions of
$\nabla^{TM}$ and $\nabla^E$ to $F$. 
The connection $\tina^F$ lifts to one on $\tilpi^*((\dol)|_F)$,
which we still denote by $\tina^F$.

We choose a local orthonormal frame such that $e_1,\dots,e_{2n_F}$
form a basis of $TF$, and $e_{2n_F+1},\cdots,e_{2n}$, that of $N$.
Denote the horizontal lift of $e_i$ ($1\le i\le 2n_F$) to $T^HN$ by $e^H_i$.
As in [\ref{BL}, Definition 8.16], we define
	\be
D^H=\sum_{i=1}^{2n_F}c(e_i)\tina^F_{e^H_i}, \quad
D^N=\sum_{i=2n_F+1}^{2n}c(e_i)\tina^F_{e_i}.
	\ee
Clearly, $D^N$ acts along the fibers of $N$.
Let $\pdrb^N$ be the $\pdrb$-operator along the fibers of $N$,
and let $(\pdrb^{N})^*$ be its formal adjoint under (\ref{inner}).
It is easy to see that $D^N=\sqrt{2}(\pdrb^N+(\pdrb^{N})^*)$.
Both $D^N$ and $D^H$ are self-adjoint with respect to (\ref{inner}).

For $T>0$, we define a scaling 
$f\in\bfh(B_{\eps_0})\mapsto S_Tf\in\bfh(B_{\eps_0\sqrt{T}})$ by
	\be
S_Tf(y,Z)=f(y,\frac{Z}{\sqrt{T}}),\quad(y,Z)\in B_{\eps_0T}.
	\ee
For a first order differential operator
	\be
Q_T=\sum_{i=1}^{2n_F}a^i_T(y,Z)\tina^F_{e^H_i}+
\sum_{i=2n_F+1}^{2n}b^i_T(y,Z)\tina^F_{e_i}+c_T(y,Z)
	\ee
acting on $\bfh(B_{\eps_0\sqrt{T}})$, where $a^i_T$, $b^i_T$ and $c_T$ are
endomorphisms of $\tilpi^*((\dol)|_F)$, we write
	\be
Q_T=O(|Z|^2\pdr^N+|Z|\pdr^H+|Z|+|Z|^p),\quad(p\in\NN)
	\ee
if there is a constant $C>0$ such that for any $T\ge1$,
$(y,Z)\in B_{\eps_0\sqrt{T}}$, we have
	\bea
\vc |a^i_T(y,Z)|\le C|Z|\quad(1\le i\le 2n_F),		\nno
\vc |b^i_T(y,Z)|\le C|Z|^2\quad(2n_F+1\le i\le 2n),	\nno
\vc |c_T(y,Z)|\le C(|Z|+|Z|^p).
	\eea

Let $J_V$ be the representation of $\mathrm{Lie}(S^1)$ on $N$.
Then $Z\mapsto J_VZ$ is a Killing vector field on $N$.
We have the following analog of [\ref{BL}, Theorem 8.18].

\begin{prop}\label{TAYLOR}
As $T\to\infty$,
	\bea
\vc S_T\,k^\hf(D^M+\ii\,Tc(V))\,k^{-\hf}S^{-1}_T			\nno
\eq \sqrt{T}(D^N+\ii\,c(J_VZ))+D^H+
	\inv{\sqrt{T}}O(|Z|^2\pdr^N+|Z|\pdr^H+|Z|+|Z|^3).
	\eea
\end{prop}

\proof{Following the proof of [\ref{BL}, Theorem 8.18], we get
	\be
S_T\,k^\hf(D^M+\ii\,Tc(V))k^{-\hf}S^{-1}_T=\sqrt{T}D^N+D^H+
	\inv{\sqrt{T}}O(|Z|^2\pdr^N+|Z|\pdr^H+|Z|).
	\ee
In fact the proof is much easier here because $F$ is totally geodesic,
hence the second fundamental form in [\ref{BL}, (8.18)] vanishes.
Next, we observe that $V=0$ on $F$ and that the vector field $J_VZ$ on $N$
is the linear approximation of $V$ on $M$ near $F$.
Further, since the actions of $S^1$ on $N$ and $M$ commute with the exponential
map, $V(y,Z)$ is odd in $Z$, and hence the second order terms vanish.
Therefore
	\be
S_T\,c(V)S^{-1}_T=\inv{\sqrt{T}}c(J_VZ)+\inv{\sqrt{T^3}}O(|Z|^3).
	\ee
The result follows.}

By Proposition~\ref{FLAT}, the solution space of the operator 
$D^N+\ii\,c(J_VZ)$ along the fiber $N_y$ ($y\in F$) is 
(the $L^2$-completion of) $K^-(N_y)\otimes E_y$.
These form an (infinite dimensional) Hermitian holomorphic bundle
$K^-(N)\otimes E|_F$ over $F$, with the Hermitian connection $\nabla^F$
induced from those in $N, E|_F\to F$.
Let $\pdrb^F$ be the corresponding operator of the twisted Dolbeault complex
$\Om^{0,*}(F,K^-(N)\otimes E|_F)$.
Set $D^F=\sum_{i=1}^{2n_F}c(e_i)\nabla^F_{e_i}=\sqrt{2}(\pdrb^F+(\pdrb^F)^*)$.

Let $\bfh^0(F)$ be the Hilbert space of square-integrable sections of
$\medwedge(T^{*(0,1)}F)\otimes K^-(N)\otimes E|_F$, and
$\bfh^0(N)$, that of the bundle $\tilpi^*((\dol)|_F)$, equipped with the
Hermitian form (\ref{inner}).
We define an embedding $\psi\colon\bfh^0(F)\to\bfh^0(N)$ by
	\be
\psi\colon\al\otimes\beta\in\bfh^0(F)\longmapsto
\tilpi^*\al\wedge\tau(\beta)\in\bfh^0(N).
	\ee
Here $\al\in\Om^{0,*}(F)$, $\beta\in L^2(K^-(N)\otimes E|_F)$
and $\tau$ is the isometry from $L^2(K^-(N)\otimes E|_F)$ to
$L^2(\tilpi^*(\medwedge(N^{*(0,1)})\otimes E|_F))$ given by
Proposition~\ref{FLAT}.
Let the image of $\psi$ be $\bfh'^{,0}=\psi(\bfh^0(F))\subset\bfh^0(N)$.
Clearly, $\psi$ is an isometry onto $\bfh'^{,0}$.
Let $p\colon\bfh^0(N)\to\bfh'^{,0}$ be the orthogonal projection.
Then we have the following analog of [\ref{BL}, Theorem 8.21].

\begin{prop}\label{REMAIN}
	\be
\psi^{-1}p\,D^Hp\,\psi=D^F.
	\ee
\end{prop}

\proof{For simplicity, we prove in the case when $E$ is a trivial line bundle
and $N$ is of rank $1$ with isotropy weight $\lam\in\ZZ\m\{0\}$.
By the definitions of $D^H$ and $D^F$, we need to prove
	\be
\tina^F_{X^H}(\tau(\beta))=\tau(\nabla^F_X\beta)
	\ee
for any $\beta\in L^2(K^-(N))$ and any vector field $X$ on $F$.
On a small neighborhood $U\subset F$, choose a unitary trivialization
$N|_U\cong U\times\CO=\set{(y,z)}{y\in U,z\in\CO}$.
Let $A\in\Om^1(U)$ be the connection $1$-form of the Hermitian connection
$\nabla^N$ in $N\to F$.
The horizontal lift of $X$ is 
$X^H=X-i_XA\,(z\frac{\pdr}{\pdr z}-\zb\frac{\pdr}{\pdr\zb})$.
If $\lam<0$, a straightforward calculation shows that for
$\beta(y,z)=f(y)\vph_k(z)$, we have
	\be
\tina^F_{X^H}\tau(\beta)(y,z)=X^H(f(y)\vph_k(z))=(Xf-k\,i_XA)(y)\vph_k(z);
	\ee
this is the connection on $(N^*)^{\otimes k}$.
If $\lam>0$, for $\beta(y,z)=f(y)\vph_\kb(z)$, we have
	\be
\tina^F_{X^H}\tau(\beta)(y,z)=(\nabla^N_X(f(y)\dr\zb)+k\,i_XA(y))\vph_\kb(z)
=(Xf+(k+1)\,i_XA)(y)\vph_\kb(z);
	\ee
this is the connection on $N^{\otimes(k+1)}$.
The result is proved.}

\subsect{3.}{Estimates of the operator and resolvent as 
$\mbox{\boldmath$\mathit T\to\infty$}$}
For $p\ge0$, let $\bfh^p(M)$, $\bfh^p(N)$ and $\bfh^p(F)$ be 
the $p$-th Sobolev spaces of the sections of the bundles $\dol\to M$,
$\tilpi^*((\dol)|_F)\to N$ and 
$\medwedge(T^{*(0,1)}F)\otimes K^-(N)\otimes E|_F\to F$, respectively.
The circle group $S^1$ acts on all these spaces; let $\bfh^p_\xi(M)$,
$\bfh^p_\xi(N)$ and $\bfh^p_\xi(F)$ be the corresponding weight spaces
of weight $\xi\in\ZZ$.
Recall the constant $\eps_0>0$ defined in the previous subsection.
We now take $\eps\in(0,\frac{\eps_0}{2}]$, which is small enough for each
eigenvalue of $L_V$ we will consider, but otherwise can be assumed fixed.
Let $\rho\colon\RE\to[0,1)$ be a smooth function such that
	\be
\rho(a)=\two{1}{a\le\hf}{0}{a\ge 1.}
	\ee
For $Z\in N$, set $\rho_\eps(Z)=\rho(\frac{|Z|}{\eps})$.
For $\beta\in L^2(K^-(N)\otimes E|_F)$, denote the component of weight 
$\xi\in\ZZ$ by $\beta_\xi\in L^2(K^-(N)\otimes E|_F)_\xi$.
Let $\al\in L^2(\medwedge(T^{*(0,1)}F))$.
We define a linear map $I\txi$ by
	\be
\sig=\al\otimes\beta\in\bfh^0(F)\longmapsto I\txi\sig=
\frac{\rho_\eps\norm{\beta}_0}{\norm{\rho_\eps S^{-1}_T(\tau(\beta_\xi))}_0}
\,\al\wedge S^{-1}_T(\tau(\beta_\xi)).
	\ee
Let the image of $I\txi$ from $\bfh^p(F)$ be
$\bfh^p\txi(N)=I\txi\bfh^p(F)\subset\bfh^p_\xi(N)$.
Denote the orthogonal complement of $\bfh^0\txi(N)$ in $\bfh^0_\xi(N)$
by $\bfh^{0,\perp}\txi(N)$,
and let $\bfh^{p,\perp}\txi(N)=\bfh^p_\xi(N)\cap\bfh^{0,\perp}\txi(N)$.
Let $p\txi$ and $p^\perp\txi$ be the orthogonal projections from
$\bfh^0_\xi(N)$ onto $\bfh^0\txi(N)$ and $\bfh^{0,\perp}\txi(N)$, respectively.

Moreover, since the bundle $\dol$ over $U_{\eps_0}$ is identified with 
$\tilpi^*((\dol)|_F)$ over $B_{\eps_0}$, we can consider $k^{-\hf}I\txi\sig$
as an element of $\bfh^p_\xi(M)$ for $\sig\in\bfh^p(F)$.
Define the linear map $J\txi$ by
	\be
\sig\in\bfh^p(F)\longmapsto J\txi\sig=k^{-\hf}I\txi\sig\in\bfh^p(M).
	\ee
Let $\bfh^p\txi(M)=J\txi\bfh^p(F)$ be the image.
Denote the orthogonal complement of $\bfh^0\txi(M)$ in $\bfh^0_\xi(M)$
by $\bfh^{0,\perp}\txi(M)$, 
and let $\bfh^{p,\perp}\txi(M)=\bfh^p_\xi(M)\cap\bfh^{0,\perp}\txi(M)$.
Let $\bar{p}\txi$ and $\bar{p}^\perp\txi$ be the orthogonal projections from
$\bfh^0_\xi(M)$ onto $\bfh^0\txi(M)$ and $\bfh^{0,\perp}\txi(M)$, respectively.
It is clear that $\bar{p}\txi=k^{-\hf}p\txi k^\hf$.

For any (possibly unbounded) operator $A$ on $\bfh^0_\xi(M)$, write
	\be
A=\left(	\begin{array}{cc} 	A^{(1)} & A^{(2)} \\
					A^{(3)} & A^{(4)}
		\end{array}	\right)
	\ee
according to the decomposition
$\bfh^0_\xi(M)=\bfh^0\txi(M)\oplus\bfh^{0,\perp}\txi(M)$, i.e.,
$A^{(1)}=\bar{p}\txi A\,\bar{p}\txi$,
$A^{(2)}=\bar{p}\txi A\,\bar{p}^\perp\txi$
$A^{(3)}=\bar{p}^\perp\txi A\,\bar{p}\txi$, and
$A^{(4)}=\bar{p}^\perp\txi A\,\bar{p}^\perp\txi$.

Let $D_T=D^M+\ii Tc(V)$.
Let $D\txi$ and $D^F_\xi$ be the restrictions of the operators $D_T$ and
$D^F$ on the subspaces $\bfh^0_\xi(M)$ and $\bfh^0_\xi(F)$, respectively,
of weight $\xi\in\ZZ$.

\begin{prop}\label{ESTIMATE}
1. As $T\to\infty$,
	\be
J^{-1}\txi D^{(1)}\txi J\txi=D^F_\xi+O\left(\inv{\sqrt{T}}\right),
	\ee
where $O(\inv{\sqrt{T}})$ denotes a first order differential operator whose
coefficients are dominated by $\frac{C}{\sqrt{T}}$ ($C>0$).\\
2. For each $\xi\in\ZZ$, there exists $C>0$ such that for any $T\ge 1$,
$\sig\in\bfh^{1,\perp}\txi(M)$, $\sig'\in\bfh^1\txi(M)$, we have
	\bea
\vc\norm{D^{(2)}\txi\sig}_0
\le C\left(\frac{\norm{\sig}_1}{\sqrt{T}}+\norm{\sig}_0\right),		\\
\vc\norm{D^{(3)}\txi\sig'}_0
\le C\left(\frac{\norm{\sig'}_1}{\sqrt{T}}+\norm{\sig}_0\right).
	\eea
3. For each $\xi\in\ZZ$, there exist $\eps\in(0,\frac{\eps_0}{2}]$, $T_0>0$,
$C>0$ such that for any $T\ge T_0$, $\sig\in\bfh^{1,\perp}\txi(M)$, we have
	\be\label{A4}
\norm{D^{(4)}\txi\sig}_0\ge C(\norm{\sig}_1+\sqrt{T}\norm{\sig}_0).
	\ee
\end{prop}

\proof{We proceed as in [\ref{BL}, Section 9].
In fact, parts 1 and 2 follow from Propositions~\ref{FLAT}, \ref{TAYLOR}
and \ref{REMAIN} the same way as in the proofs of [\ref{BL}, Theorems 9.8, 9.9]
with trivial modifications.
Part 3 is the analog of [\ref{BL}, Theorem 9.14].
The key observation is that we are restricting the operators to a fixed weight
space $\bfh_\xi(M)$, where $L_V$ is a constant.
Thus the terms
	\be\label{bounded}
-2\ii L_V+\hf\ii\sumi c(e_i)c(\nabla_{e_i}V)
+\ii\,\tr\,\nabla.V|_{T\ah M}-2\ii\,r_V
	\ee
in (\ref{laplace}) is bounded on $\bfh_\xi(M)$.
This is the analog of the fact that $[D,V]$ is of order zero in [\ref{BL}].
The boundedness of (\ref{bounded}) allows us to localize the problem
to sufficiently small neighborhood $U_\eps$ of $F$ for each weight $\xi\in\ZZ$.
Using Propositions~\ref{FLAT}, \ref{TAYLOR} and \ref{REMAIN} again, 
(\ref{A4}) can be proved by similar arguments as in [\ref{BL}, Section 9b).3].
The details are left to the interested reader.}

For two bounded operators $A\in{\cal L}(\bfh^0_\xi(M))$,
$B\in{\cal L}(\bfh^0_\xi(F))$, set
	\be
\dr(A,B)=\sum_{j=2}^4\norm{A^{(j)}}_1+\norm{J^{-1}\txi A^{(1)}J\txi-B}_1.
	\ee
We fix a constant $c_0\in(0,1]$ such that
	\be
\spec(D^F_\xi)\cap[-2c_0,2c_0]\subset\{0\},
	\ee
where $\spec$ denotes the spectrum of an operator.
Then we have the following analog of [\ref{BL}, Theorem 9.23].

\begin{prop}\label{RESOLVANT}
For any $\xi\in\ZZ$, there exists $T_0\ge 1$, such that for any $T\ge T_0$,
$\lam\in\CO$ with $|\lam|=c_0$, $\lam-D\txi$ is invertible on $\bfh^0(M)$.
Moreover, for any integer $p\ge 2n+2$, there exists $c_p>0$ such that for
any $T\ge T_0$, $\lam\in\CO$ with $|\lam|=c_0$, we have
	\be
\dr((\lam-D\txi)^{-p},(\lam-D^F_\xi)^{-p})\le\frac{c_p}{\sqrt{T}}.
	\ee
\end{prop}

\proof{In view of Proposition~\ref{ESTIMATE}, the result can be proved using
the same formal arguments in [\ref{BL}, Sections 9c)-9e)].}

\subsect{4.}{Proof of Theorem~\ref{CIRCLE}}
Let $\gam$ be the circle in $\CO$ of center $0$ and radius $c_0$, oriented
counterclockwise when necessary.
By Proposition~\ref{RESOLVANT}, $\gam\cap\spec(D\txi)=\emptyset$ for $T$ large
enough.
Let $K^{c_0}\txi$ be the direct sum of eigenspaces of $D\txi$ associated to
the eigenvalues $\lam$ such that $|\lam|\le c_0$.
For $T$ large enough,
	\be\label{proj}
P^{c_0}\txi=\inv{2\pi\ii}\int_\gam(\lam-D\txi)^{-1}\dr\lam
	\ee
is the orthogonal projection from $\bfh^0_\xi(M)$ onto $K^{c_0}\txi$.
Integrating by parts in (\ref{proj}), we get for any $p\in\NN$,
	\be
P^{c_0}\txi=\inv{2\pi\ii}\int_\gam\lam^{p-1}(\lam-D\txi)^{-p}\dr\lam.
	\ee
Using Proposition~\ref{RESOLVANT}, we obtain, for some $C>0$,
	\be
\dr(P^{c_0}\txi,P^F_\xi)\le\frac{C}{\sqrt{T}},
	\ee
where $P^F_\xi$ is the orthogonal projection from $\bfh^0_\xi(F)$ onto 
$K_\xi=\ker D^F_\xi$.
Therefore for $T$ large enough,
	\be\label{ker}
\dim_\co K^{c_0}\txi=\dim_\co K_\xi=
\sum_{k=0}^{n_F}\dim_\co H^k(F,\OO(K^-(N)\otimes E|_F)_\xi).
	\ee
In fact, (\ref{ker}) holds for each degree in the Dolbeault complexes.
By Hodge theorem and Proposition~\ref{FLAT}, we have for any $0\le k\le n_F$,
	\be
J\txi\,H^k(F,\OO(K^-(N)\otimes E|_F)_\xi)\subset
\bfh^0_\xi(M)\cap L^2(\medwedge^{k+(n-n_F-\nu_F)}(T^{*(0,1)}M)\otimes E).
	\ee
Let $m^k_\xi$ be the sum of the multiplicities of all eigenvalues $\lam$ of
$D\txi$ on $(0,k)$-forms such that $|\lam|<c_0$.
Then
	\be\label{small}
m^k_\xi=\dim_\co H^{k-(n-n_F-\nu_F)}(F,\OO(K^-(N)\otimes E|_F)_\xi).
	\ee
On the other hand, an obvious application of (the equivariant version of)
the Hodge theorem yields that $m^k_\xi$ ($0\le k\le n$) satisfy
the Morse-type inequalities
	\be\label{hodge-morse}
\sumk t^km^k_\xi=\sumk t^k\mult_\xi H^k(M,\OO(E))+(1+t)Q_\xi(t),
	\ee
where $Q_\xi(t)\ge 0$.
Combining (\ref{small}) and (\ref{hodge-morse}), we get
	\be
\sumk t^k\mult_\xi H^{k-(n-n_F-\nu_F)}(F,\OO(K^-(N)\otimes E|_F))
=\sumk t^k\mult_\xi H^k(M,\OO(E))+(1+t)Q_\xi(t).
	\ee
This is exactly (\ref{strong1-}).	$\hfill\Box$

\sect{Applications}
We consider various applications of the equivariant holomorphic Morse
inequalities.
In subsection 4.1, we apply the inequalities (\ref{strong})
to exterior powers of holomorphic cotangent bundles.
The results are relations of the Hodge numbers of compact \ka manifolds
and the fixed submanifolds of torus actions.
In subsection 4.2, we obtain a gluing formula of the Poincar\'e-Hodge
polynomials in the context of symplectic cutting.
We also study geometric quantization by applying the inequalities to
the prequantum line bundles.

\subsect{1.}{Relations among Hodge numbers of the fixed submanifolds}
Consider a compact \ka manifold $(M,\om)$ of complex dimension $n$.
Let $H\kl(M)=H^l(M,\EE^k(M))$ ($k,l=0,1,\dots,n$) be the Dolbeault
cohomology groups of $M$, where $\EE^k(M)=\OO(M,\medwedge^k(T^{*(1,0)}M))$
is the sheaf of holomorphic sections in the $k$-th exterior power of 
the holomorphic cotangent bundle of $M$.
Let $h\kl(M)=\dim_\co H\kl(M)$ be the {\em Hodge numbers} of $M$ and
$P(M;s,t)=\sum^n_{k,l=0}s^kt^lh\kl(M)$, the {\em Poincar\'e-Hodge polynomial}.
(The Poincar\'e polynomial of $M$ is $P(M;t)=P(M;t,t)$.)
The {\em Todd genus} (or {\em arithmetic genus}) of $M$ is given by 
$\tau(M)=P(M;0,-1)$.
We also have the well-known relations $h\kl(M)=h^{l,k}(M)=h^{n-k,n-l}(M)$
for any $0\le k,l\le n$, or $P(M;s,t)=P(M;t,s)=(st)^nP(M;s^{-1},t^{-1})$
for compact \ka manifolds.

Suppose that there is a holomorphic action of the torus group $T$ on $M$
preserving the \ka form.
As before, if the fixed-point set $F$ is non-empty, it is a disjoint union
of connected compact \ka submanifolds $F_1,\dots,F_m$, of (complex) 
dimensions $n_1,\dots,n_m$, respectively.
Recall that the weights $\lam_{r,k}$ ($1\le r\le m$, $1\le k\le n$)
of isotropy representations of $T$ on the normal bundles of $F_r$
cut the Lie algebra $\gt$ into action chambers.

\begin{thm}\label{HODGE}
1. All the cohomology groups $H\kl(M)$ ($k,l=0,1,\dots,n$) are trivial
representations of $T$.\\
\noindent 2. For any choice of action chamber $C$,
	\be\label{hodge}
\sumr(st)^{\ncr}P(F_r;s,t)=P(M;s,t).
	\ee
In particular, if all the fixed points are isolated, then for any choice of
chamber $C$, we have
	\be\label{isohodge}
h^{k,k}(M)=\#\set{p\in F}{\nu^C_p=k} \quad (0\le k\le n)
	\ee
and $h\kl(M)=0$ for $k\ne l$.
\end{thm}

\proof{The $T$-action can be lifted holomorphically to the bundles
$\medwedge^k(T^{*(1,0)}M)$ ($0\le k\le n$).
Applying (\ref{strong}), we obtain, for some $Q^{k,C}(t)\ge 0$,
	\be\label{ext}
\sumr t^\ncr\sum^{n_r}_{l=0}t^l\ch H^l(F_r,\OO(T^{k,C}_r))
=\sum^n_{l=0}t^l\ch H\kl(M)+(1+t)Q^{k,C}(t),
	\ee
where
	\be
T^{k,C}_r=S((\Ncr)^*)\otimes S(\Nmcr)\otimes\medwedge^{\ncr}(\Nmcr)
\otimes\medwedge^k(N^{*(1,0)}_r\oplus T^{*(1,0)}F_r).
	\ee
If $\lam_{r,l}$ ($1\le l\le n-n_r$) are the weights of the $T$-action on
the fiber of $N^{1,0}$, then the weights in $T^k_r$ are of the form
	\be\label{weight}
-\sum^{n-n_r}_{l=1}m_l\lam^C_{r,l}+\sum_{\lam_{r,l}\in-C^*}\lam_{r,l}
-\sum_{l\in I}\lam_{r,l},
	\ee
where $m_l\ge0$ and $I\subset\{1,\dots,n\}$ with $|I|\le k$.
(\ref{weight}) is an element in the closed cone $-\overline{C^*}$;
it is $0$ if and only if all $m_l=0$ and $I=\set{l}{\lam_{r,l}\in-C^*}$ 
(hence $|I|=\ncr\le k$), corresponding to the sub-bundle 
$\medwedge^{k-\ncr}(T^{*(1,0)}F_r)$ (if $k\ge\ncr$) of $T^{k,C}_r$.
By the weak inequalities (\ref{weak}), we conclude that 
$\supp H\kl(M)\subset-\overline{C^*}\cap\lat^*$.
Had we chosen the opposite chamber $-C$, we would get 
$\supp H\kl(M)\subset\overline{C^*}\cap\lat^*$.
Therefore $\supp H\kl(M)\subset\{0\}$, and hence $H\kl(M)$ are trivial 
representations of $T$.
Restricting (\ref{ext}) to $0\in\lat^*$, we get
	\be
\sumr t^\ncr\sum^{n_r}_{l=0}t^lh^{k-\ncr,l}(F_r)
=\sum^n_{l=0}t^lh\kl(M)+(1+t)Q^{k,C}_0(t)
	\ee
for some polynomial $Q^{k,C}_0(t)\ge0$ in $\ZZ[t]$.
By the symmetry of Hodge numbers, we obtain
	\be\label{ineq}
\sumr(st)^{\ncr}P(F_r;s,t)=P(M;s,t)+(1+s)(1+t)Q^C_0(s,t),
	\ee
where $Q^C_0(s,t)\ge0$ in $\ZZ[s,t]$ is defined by
the relation $\sumk s^kQ^{k,C}_0(t)=(1+s)Q^C_0(s,t)$.
Now let $s=t$ in (\ref{ineq}). 
Since the moment map (after projecting to any direction in $\pm C$) is
a perfect Morse function in the sense of Bott (see for example [\ref{Fr}]),
we get $Q^C_0(t,t)=0$.
Therefore $Q^C_0(s,t)=0$ and (\ref{hodge}) follows.}

\begin{rmk} {\em
1. Part 1 of Theorem~\ref{HODGE} is in fact true for any connected
group acting holomorphically and isometrically:
Since the action of any group element is homotopic to the identity map,
its actions on the de Rham cohomology groups $H^k(M)$ ($0\le k\le 2n$)
are trivial.
Moreover, since the action preserves the complex and \ka structures,
it preserves the Hodge decomposition. Hence the result.\\
2. Part 2 of the theorem in the case of $\CO^*$-actions was obtained by
Carrell-Sommese [\ref{CS}, Theorem 2].
For torus actions, that formula (\ref{hodge}) is true for any choice of
action chamber gives various constraints on the fixed-point data.
For example, when the fixed-points are isolated, the number of fixed points
$p\in F$ with polarizing index $\nu^C_p=k$ is independent of the choice of 
the chamber $C$.
In general, consider two action chambers $\pm C$. We have
	\be\label{consistent}
\sumr(st)^{\nu^{-C}_r}P(F_r;s,t)=\sumr(st)^{\ncr}P(F_r;s,t).
	\ee
In fact, using $\nu^{-C}_r=n-n_r-\ncr$ and 
$P(F_r;s,t)=(st)^{n_r}P(F_r;s^{-1},t^{-1})$, it is not hard to see
that (\ref{consistent}) is equivalent to $P(M;s,t)=(st)^nP(M;s^{-1},t^{-1})$.}
\end{rmk}

\begin{cor}\label{todd}
For any choice of action chamber $C$, there is only one component $F_{r_C}$
with $\nu^C_{r_C}=0$.
Moreover, $h^{0,k}(M)=h^{0,k}(F_{r_C})$ for all $0\le k\le n$.
In particular, $\tau(M)=\tau(F_{r_C})$.
\end{cor}

\proof{The first part is well-known [\ref{A},\ref{GS}] and is included
here for completeness;
in fact, $F_{r_C}$ is the inverse image of a vertex of the moment polytope.
Taking $s=0$ in (\ref{hodge}), we get $P(M;0,t)=P(F_{r_C};0,t)$.
The rest follows easily.}

\begin{ex} {\em 
Hamiltonian $S^1$-actions on symplectic $4$-manifolds have been classified
up to $S^1$-equivariant diffeomorphisms [\ref{Au}], and subsequently,
up to $S^1$-equivariant symplectomorphisms [\ref{Kar}].
Moreover it was shown that all such manifolds are K\"ahler [\ref{Kar}].
Let $M$ be such a manifold with a Hamiltonian $S^1$-action and let
$h\colon M\to\RE$ be a moment map.
Suppose $\Sig_\pm$ are the critical components on which $h$ reaches its
maximum, minimum, respectively.
By dimensional reasons, other critical points of $h$ have Morse indices $2$
and are isolated; let $m_2$ be the number of such.
The Hodge-Poincar\'e polynomial of $M$ is
	\bea
P(M;s,t)\eq P(\Sig_-;s,t)+stP(\Sig_+;s,t)+m_2st		\nno
	\eq P(\Sig_+;s,t)+stP(\Sig_-;s,t)+m_2st
	\eea
It follows that there are three possibilities:\\
\noindent I: Both $\Sig_\pm$ are isolated points, i.e., $P(\Sig_\pm;s,t)=1$.
In this case, $h^{0,0}(M)=h^{2,2}(M)=1$, $h^{1,1}(M)=m_2=b_2(M)$,
others$\,=0$.\\
\noindent II: One of them is reached at an isolated point, the other at 
a sphere, i.e., $P(\Sig_\pm;s,t)=1$ and $st$, respectively. 
In this case, $h^{0,0}(M)=h^{2,2}(M)=1$, $h^{1,1}(M)=m_2+1=b_2(M)$,
others$\,=0$.\\
\noindent III: Both $\Sig_\pm$ are Riemann surfaces of the same genus, 
i.e., $P(\Sig_\pm;s,t)=1+g(s+t)+st$, where $g$ is the genus.
In this case, $h^{0,0}(M)=h^{2,2}(M)=1$,
$h^{0,1}(M)=h^{1,0}(M)=h^{1,3}(M)=h^{3,1}(M)=g=\hf b_1(M)$,
$h^{1,1}(M)=m_2+2=b_2(M)$, others$\,=0$.\\
\noindent In all cases, $b^+_2(M)=1$.}
\end{ex}

If $V$ is a holomorphic Killing vector field on $(M,\om)$,
let $Z_1,\dots,Z_m$ be the connected components of the zero-set of $V$,
of (complex) dimensions $n_1,\dots,n_m$, respectively.
We have the following special case of the theorem of Carrell-Lieberman 
[\ref{CL},\ref{CS}], which holds for a more general (not necessarily Killing)
holomorphic vector field.

\begin{cor}
$h\kl(M)=0$ if $|k-l|>\max\{n_1,\dots,n_m\}$.
\end{cor}

\proof{The $1$-parameter group generated by $V$ is a subgroup
of the isometry group of $M$.
Therefore its closure is a torus $T$, whose fixed-point set is precisely $Z$.
We choose any action chamber $C$ (for example, the one which contains the
generator $V$).
Let $s=t^{-1}$ in (\ref{hodge}). We get
	\be
\sumr\sum^{n_r}_{k,l=0}h\kl(Z_r)t^{l-k}=\sum^n_{k,l=0}h\kl(M)t^{l-k}.
	\ee
The result follows easily.}

When the fixed-point set $F$ is discrete, this result was deduced
from holomorphic Morse inequalities in [\ref{W}].

\subsect{2.}{Symplectic quotients, symplectic cuts and quantization}
Let $(M,\om)$ be a symplectic manifold of (real) dimension $2n$.
If the circle group $S^1$ acts Hamiltonianly on $(M,\om)$,
let $V$ be the vector field on $M$ that generates the action
and $h\colon M\to\RE$, the moment map such that $i_V\om=\dr h$.
If $0$ is a regular value of $h$, then $S^1$ acts locally freely on 
$h^{-1}(0)$ and the symplectic quotient $M_0=h^{-1}(0)/S^1$ is an orbifold.
Let $i\colon h^{-1}(0)\to M$ be the inclusion and 
$\pi\colon h^{-1}(0)\to M_0$, the projection.
There is a canonical symplectic form $\om_0$ on $M_0$ such that
$\pi^*\om_0=i^*\om$.
To avoid orbifold singularities, we further assume that $S^1$ acts freely on 
$h^{-1}(0)$.
In this case, $M_0$ is a smooth manifold of dimension $2n-2$.

We now recall the notion of symplectic cutting [\ref{L}].
Let the complex plane $\CO$ be equipped with the standard \ka form
$\om=\frac{\ii}{2}\dr z\wedge\dr\zb$.
Consider two actions of the circle group $S^1$ on $\CO$ with weights $\pm 1$.
Both actions are Hamiltonian with the moment maps $\mp\hf|z|^2$, respectively.
The diagonal actions of $S^1$ on $M\times\CO$ are again Hamiltonian and
the moment maps are $h\mp\hf|z|^2$, of which $0$ is still a regular value.
Let $M_\pm$ be the symplectic quotients of $M\times\CO$ at level $0$
by the two $S^1$-actions defined above.
$(M_\pm,\om_\pm)$ are smooth symplectic manifolds with Hamiltonian 
$S^1$-actions and $M_0$ is embedded as one of the components (still denoted by
$M_0$) fixed by  $S^1$; the compliments $M_\pm\m M_0$ are $S^1$-equivariantly
symplectomorphic to $h^{-1}(\RE^\pm)\subset M$, respectively.
Therefore the fixed-point set of $M_\pm$ is the union of $M_0$ and all
the components $F_r$ ($1\le r\le m$) such that $h(F_r)\in\RE^\pm$.
Moreover the circle bundle of the normal bundle $N_0$ of $M_0$ in $M_+$
is isomorphic to the bundle $h^{-1}(0)\to M_0$, with weight $-1$, while 
the normal bundle of $M_0$ in $M_-$ is isomorphic to $N^*_0$, with weight $1$.

When $(M,\om)$ is a \ka manifold with a Hamiltonian $S^1$-action,
then the symplectic quotient $(M_0,\om_0)$ and two symplectic cuts
$(M_\pm,\om_\pm)$ are also K\"ahler.
We now establish a gluing formula for Poincar\'e-Hodge polynomials.

\begin{prop}
Let $(M,\om)$ be a compact \ka manifold with a Hamiltonian $S^1$-action.
If $0$ is a regular value of the moment map $h$ and $S^1$ acts freely on
$h^{-1}(0)$, then
	\be\label{glue}
P(M_+;s,t)+P(M_-;s,t)=P(M;s,t)+(1+st)P(M_0;s,t).
	\ee
\end{prop}

\proof{If we choose the positive chamber chamber $C=\RE^+$, then $M_0$
is embedded in $M_\pm$ with polarizing indices $1$ and $0$, respectively.
Applying (\ref{hodge}) to $M_\pm$, we get
	\be
P(M_+;s,t)=\sum_{F_r\subset h^{-1}(\re^+)}(st)^{\nu_r}P(F_r;s,t)+stP(M_0;s,t)
	\ee
and
	\be
P(M_-;s,t)=\sum_{F_r\subset h^{-1}(\re^-)}(st)^{\nu_r}P(F_r;s,t)+P(M_0;s,t).
	\ee
(\ref{glue}) follows from the above two equalities and (\ref{hodge}).}

As a consequence, we have the following interesting 

\begin{cor}
Under the above assumptions, $h^{0,k}(M)=h^{0,k}(M_0)$ for all $0\le k\le n$.
Hence $\tau(M)=\tau(M_0)$.
\end{cor}

\proof{From Corollary~\ref{todd} we get $P(M_+;0,t)=P(M_-;0,t)=P(M_0;0,t)$.
This, together with (\ref{glue}) when $s=0$, implies that
$P(M;0,t)=P(M_0;0,t)$.}

By induction, these results are true for Hamiltonian torus actions.
In fact, the second part is true when $M$ is a general symplectic (hence
almost complex) manifold with a (possibly) non-Abelian group action
[\ref{MS}, \ref{TZ}].
The first part is the refinement of this result when $M$ is K\"ahler. 

Recall that a symplectic manifold $(M,\om)$ is {\em quantizable}
if the de Rham class $[\frac{\om}{2\pi}]\in H^2(M,\ZZ)$.
In this case, there is a complex line bundle, called the {\em pre-quantum
line bundle}, with a connection $\nabla$ whose curvature is $\frac{\om}{\ii}$.
We have the following

\begin{lemma}\label{PREQ}
Suppose that a symplectic manifold $(M,\om)$ is quantizable and is equipped
with a Hamiltonian $S^1$-action with moment map $h$.
If one of $h(F_r)\in\ZZ$, then all $h(F_r)\in\ZZ$ and
the $S^1$-action can be lifted to the pre-quantum line bundle $L$.
In this case, $M_\pm$, $M_0$ are quantizable; let $L_\pm$, $L_0$ 
be their pre-quantum line bundles.
There are following isomorphisms of line bundles with connections:
$L_\pm|_{M_0}\cong L_0$,
$L_\pm|_{M_\pm\m M_0}\cong L|_{h^{-1}(\re^\pm)}$.
\end{lemma}

\proof{The generator of $\mathrm{Lie}(S^1)$ acts on 
the space of sections of $L$ by $-\nabla_V+\ii h$;
this gives an $\RE$-action on $L$ preserving $\nabla$.
If $h(F_r)\in\ZZ$, then the $\RE$-action is an $S^1$-action on $L|_{F_r}$.
Since $\nabla$ is $\RE$-invariant, the parallel transport commutes with
the $\RE$-action.
Therefore the $\RE$-action factorizes through $S^1$ on the total space of $L$.
In particular, $h(F_r)\in\ZZ$ on any fixed component $F_r$.
The line bundle $L_0=i^*L/S^1\to M_0$ has a connection $\nabla^0$ such that
$\pi^*\nabla^0=i^*\nabla$ so that its curvature is $\frac{\om_0}{\ii}$
[\ref{GS2}].
Following the construction of symplectic cuts, there are $S^1$-invariant
pre-quantum line bundles $L_\pm$ on $M_\pm$ whose curvature is 
$\frac{\om_\pm}{\ii}$.
The isomorphisms in the last part were proved in [\ref{DGMW},\ref{M}].}

If in addition $(M,\om)$ is K\"ahler, the pre-quantum line bundle $L$
can be made into an $S^1$-invariant holomorphic Hermitian line bundle.
Therefore the line bundle $L_0$ over $M_0$ is also holomorphic and Hermitian,
and so are $L_\pm$.
Furthermore, the actions of $S^1$ on $L_\pm$ preserve holomorphic structures.
By {\em quantization} on $(M,\om)$ we mean to associate to $(M,\om)$
the virtual vector space
	\be\label{virtual}
\HH(M)=\bigoplus^n_{k=0}(-1)^k H^k(M,\OO(L)).
	\ee
When $M$ is not a complex manifold, each individual cohomology group in
(\ref{virtual}) does not make sense, however the alternating sum can be
defined as the index of a spin$^\co$-Dirac operator using only an almost
complex structure.
In [\ref{DGMW}], it was proved that under this more general setting, we have
the following relation on quantization, symplectic cutting and reduction
	\be\label{com}
\dim_\co\HH(M_\pm)^{S^1}=\dim_\co\HH(M_0)=\dim_\co\HH(M)^{S^1},
	\ee
and the gluing formula
	\be\label{dgmw}
\ch\HH(M)=\ch\HH(M_+)+\ch\HH(M_-)-\dim_\co\HH(M_0).
	\ee

The last equality in (\ref{com}) was the $S^1$-case of a conjecture by
Guillemin and Sternberg [\ref{GS2}]; the case with compact non-Abelian
group actions was proved by Meinrenken [\ref{M}], and Vergne [\ref{V}],
Jeffrey and Kirwan [\ref{JK}] and others under various generalities
using localization techniques,
and by Tian and Zhang [\ref{TZ}] using a direct analytic approach.
Moreover, when $M$ is a compact \ka manifold with a Hamiltonian action of
(possibly) non-Abelian group $G$, there are Morse-type inequalities 
which bound the invariant cohomologies of $M$ in terms of those of 
the symplectic quotient $M_0$ (of complex dimension $n_0$)[\ref{TZ}].
That is,
	\be\label{tz}
\sum_{k=0}^{n_0}t^k\dim_\co H^k(M_0,\OO(L_0))=
\sumk t^k\dim_\co H^k(M,\OO(L))^G+(1+t)Q_0(t)
	\ee
for some $Q_0(t)\ge0$.
We prove similar Morse-type inequalities relating quantizations on 
$M_0$ and $M_\pm$, which will strengthen the first equality of (\ref{com})
when $M$ is K\"ahler and the group is $S^1$.

\begin{prop}
If $M$ is K\"ahler, then there exist polynomials $Q^\pm_0(t)\ge0$ such that
	\be\label{pm}
\sum_{k=0}^{n-1}t^k\dim_\co H^k(M_0,\OO(L_0))=
\sumk t^k\dim_\co H^k(M_\pm,\OO(L_\pm))^{S^1}+(1+t)Q^\pm_0(t).
	\ee
\end{prop}

\proof{Consider the $S^1$-action on $M_+$,
whose fixed-point set consists of $M_0$ and $F_r$ with $h(F_r)>0$.
The weights of the $S^1$-action on the fibers of $L_+$ over
$M_0$ and $F_r$ are $0$ and $h(F_r)$, respectively.
Applying (\ref{strong1-}) to $M_+$, we obtain
	\bea\label{pol}
\vc\sum_{k=0}^{n-1}t^k\ch H^k(M_0,\OO(S(N^*_0)\otimes L_0))
  +\sum_{F_r\subset h^{-1}(\re^+)}t^k\ch H^k(F_r,\OO(K^-(N_r)\otimes L|_{F_r}))
						\nno
\vc\quad=\sumk t^k\ch H^k(M_+,\OO(L_+))+(1+t)Q^+(t),
	\eea
for $Q^+(t)\ge 0$.
Notice that all the weights on $S(N^*_0)\otimes L_0$ and 
$K^-(N_r)\otimes L|_{F_r}$ are non-negative,
and the zero weight comes only from the sub-bundle $L_0$ of the former.
By restricting (\ref{pol}) to the zero weight, we obtain (\ref{pm}) for $M_+$.}

\begin{rmk} {\em In the light of (\ref{dgmw}), we conjecture that
there is an $S^1$-equivariant Mayer-Vietoris-type long exact sequence
	\be
\cdots\to H^k(M,\OO(L))\to H^k(M_+,\OO(L_+))\oplus H^k(M_-,\OO(L_-))\to
H^k(M_0,\OO(L_0))\to H^{k+1}(M,\OO(L))\to\cdots
	\ee
when $M$ is K\"ahler.
If so, then there is a polynomial $Q(t)\ge0$ such that
	\bea\label{exact}
\vc\quad\sumk t^k\ch H^k(M,\OO(L))
	+\sum_{k=0}^{n-1}t^k\dim_\co H^k(M_0,\OO(L_0))		\nno
\eq\sumk t^k\ch H^k(M_+,\OO(L_+))+\sumk t^k\ch H^k(M_-,\OO(L_-))+(1+t)Q(t).
	\eea
In fact, the polynomial $Q(t)\le\sum_{k=0}^{n-1} t^k\ch H^k(M_0,\OO(L_0))$.
Therefore (\ref{exact}), if correct, also implies
	\bea
\vc\sumk t^k\ch H^k(M_+,\OO(L_+))+\sumk t^k\ch H^k(M_-,\OO(L_-))
   +\sum_{k=0}^{n-1} t^{k+1}\ch H^k(M_0,\OO(L_0))			\nno
\eq\quad\sumk t^k\ch H^k(M,\OO(L))+(1+t)Q'(t)
	\eea
for some $Q'(t)\ge0$.}
\end{rmk}

\noindent {\bf Acknowledgement.} 
Part of the work of S.W.\ was done in MSRI, Berkeley,
supported by NSF grant DMS-90-22140.
S.W.\ thanks the hospitality of Nankai Institute of Mathematics
where this work was completed.
Part of the work of W.Z.\ was done in Courant Institute, New York.
W.Z. is supported in part by the Chinese National Science Foundation
and the State Education Commission of China.

\bigskip

        \newcommand{\athr}[2]{{#1}.\ {#2}}
        \newcommand{\au}[2]{\athr{{#1}}{{#2}},}
        \newcommand{\an}[2]{\athr{{#1}}{{#2}} and}
        \newcommand{\jr}[6]{{#1}, {\it {#2}} {#3}\ ({#4}) {#5}-{#6}}
        \newcommand{\pr}[3]{{#1}, {#2} ({#3})}
        \newcommand{\bk}[4]{{\it {#1}}, {#2}, ({#3}, {#4})}
        \newcommand{\cf}[8]{{\it {#1}}, {#2}, {#5},
                 {#6}, ({#7}, {#8}), pp.\ {#3}-{#4}}
        \vspace{5ex}	
	\newpage
        \begin{flushleft}
{\bf References}
        \end{flushleft}
{\small
        \begin{enumerate}
        
        \item\label{A}
        \au{M.\ F}{Atiyah}
        \jr{Convexity and commuting Hamiltonians}
        {Bull. London Math. Soc.}{14}{1982}{1}{15}

        \item\label{AB}
        \an{M.\ F}{Atiyah} \au{R}{Bott}
        \jr{A Lefschetz fixed point formula for elliptic complexes, Part I}
        {Ann. Math.}{86}{1967}{374}{407};
        \jr{Part II}{Ann. Math.}{87}{1968}{451}{491}

	\item\label{ASS}
	\an{M.\ F}{Atiyah} \au{G}{Segal}
	\jr{The index of elliptic operators. II}
	{Ann. Math.}{87}{1968}{531}{545};\\
	\an{M.\ F}{Atiyah} \au{I.\ M}{Singer}
	\jr{The index of elliptic operators. III}
	{Ann. Math.}{87}{1968}{546}{604}

        \item\label{Au}
        \au{M}{Audin}
        \cf{Hamiltoniens p\'eriodiques sur les vari\'et\'es symplectiques
            compactes de dimension 4}
           {G\'eom\'etrie symplectique et m\'echanique, Proceedings 1988}
           {1}{25}{C.$\,$Albert ed., Lecture Notes in Mathematics 1416}
           {Springer}{Berlin, Heidelberg, New York}{1990};
        \bk{The topology of torus action on symplectic manifolds,
           {\rm Progress in Mathematics, Vol.$\,$93}}
           {Birkh\"{a}user}{Basel, Boston, Berlin}{1991}, Ch.$\,$IV

	\item\label{Bi} 
	\au{J.-M}{Bismut} 
	\jr{The Witten complex and the degenerate Morse inequalities} 
	{J. Diff. Geom.}{23}{1986}{207}{240}

	\item\label{BL} 
	\an{J.-M}{Bismut} \au{G}{Lebeau}
	\jr{Complex immersions and Quillen metrics}
	{Inst. Hautes Etudes Sci. Publ. Math.}{74}{1991}{1}{297}

	\item\label{CG}
	\an{A}{Canas da Silva} \au{V}{Guillemin}
	\pr{On the Kostant multiplicity formula for group actions with
	non-isolated fixed points}{MIT math preprint}{1996}

	\item\label{CL}
	\an{J.\ B}{Carrell} \au{D}{Lieberman}
	\jr{Holomorphic vector fields and Kaehler manifolds}
	{Invent. Math.}{23}{1973}{303}{309}

	\item\label{CS}
	\an{J.\ B}{Carrell} \au{A.\ J}{Sommese}
	\jr{Some topological aspects of ${\bf C}^*$ actions on compact Kaehler
	manifolds}{Comment. Math. Helvetica}{54}{1979}{567}{582}

%
        \item\label{DGMW}
        \au{H}{Duistermaat} \au{V}{Guillemin} \an{E}{Meinrenken} \au{S}{Wu}
        \jr{Symplectic reduction and Riemann-Roch for circle actions}
        {Math. Res. Lett.}{2}{1995}{259}{266}

	\item\label{Fr}
	\au{T}{Frankel}
	\jr{Fixed points and torsion on \ka manifold}
	{Ann. Math.}{70}{1959}{1}{8}

 	\item\label{GLS}
	\au{V}{Guillemin} \an{E}{Lerman} \au{S}{Sternberg}
	\jr{On the Kostant multiplicity formula}
	{J. Geom. Phys.}{5}{1988}{721}{750}

	\item\label{GS}
        \an{V}{Guillemin} \au{S}{Sternberg}
        \jr{Convexity properties of the moment mapping}
        {Invent. Math.}{67}{1982}{491}{513}

        \item\label{GS2}
        \an{V}{Guillemin} \au{S}{Sternberg}
        \jr{Geometric quantization and multiplicities of group representations}
        {Invent. Math.}{67}{1982}{515}{538}

	\item\label{JK}
	\an{L.\ C}{Jeffrey} \au{F.\ C}{Kirwan}
	\jr{Localization and the quantization conjecture}
	{Topology}{36}{1997}{647}{693}

	\item\label{Kar}
	\au{Y}{Karshon}
	\pr{Periodic Hamiltonian flows on four dimensional manifolds}
	{MIT math preprint}{1994}

        \item\label{L}
        \au{E}{Lerman}
        \jr{Symplectic cuts}
        {Math. Res. Lett.}{2}{1995}{247}{258}

        \item\label{MW}
        \an{V}{Mathai} \au{S}{Wu}
        \pr{Equivariant holomorphic Morse inequalities I: a heat kernel proof}
        {ICTP preprint IC/96/29, {\tt dg-ga/9602007}}{1996}

        \item\label{M}
        \au{E}{Meinrenken}
        \pr{Symplectic surgery and the spin$^c$-Dirac operator}
        {MIT math. preprint}{1995}, to appear in Adv. Math.

        \item\label{MS}
        \an{E}{Meinrenken} \au{R}{Sjamaar}
	\pr{Singular reduction and quantization}
	{preprint}{1996}

	\item\label{PW}
	\an{E}{Prato} \au{S}{Wu}
	\jr{Duistermaat-Heckman measures in a non-compact setting}
	{Comp. Math.}{94}{1994}{113}{128}

	\item\label{TZ}
	\an{Y}{Tian} \au{W}{Zhang}
	\pr{Symplectic reduction and quantization,
	Symplectic reduction and analytic localization}
	{preprints}{1996}

	\item\label{V}
	\au{M}{Vergne}
	\jr{Multiplicities formula for geometric quantization. I, II}
	{Duke Math. J.}{82}{1996}{143}{179}, 181-194

	\item\label{W0}
	\au{E}{Witten} 
	\jr{Supersymmetry and Morse theory}
	{J. Diff. Geom.}{17}{1982}{661}{692}
	
	\item\label{W} 
	\au{E}{Witten} 
	\cf{Holomorphic Morse inequalities}
	{Algebraic and differential topology, Teubner-Texte Math., 70}
	{318}{333}{ed.\ G.\ Rassias}{Teubner}{Leipzig}{1984}

	\item\label{WuII}
	\au{S}{Wu}
	\pr{Equivariant holomorphic Morse inequalities II: torus and
	non-Abelian group actions}
	{MSRI preprint No.~1996-013, {\tt dg-ga/9602008}}{1996}

        \end{enumerate}}
	\end{document}